\documentclass{svmult-ddm}

\usepackage[utf8]{inputenc}
\usepackage{helvet}         
\usepackage{courier}        
\usepackage{graphicx}        

\usepackage[bottom]{footmisc}

\usepackage{amssymb}
\usepackage{amsfonts}
\usepackage{amsbsy}
\usepackage{amscd}
\usepackage{amstext}
\usepackage{dsfont}
\usepackage[english]{babel}
\usepackage{color}
\usepackage{graphics}
\usepackage{epsfig}
\usepackage{subfigure}
\usepackage{wrapfig}
\usepackage{psfrag}
\usepackage{color}
\usepackage{url}
\usepackage{verbatim}
\usepackage{algorithm}
\usepackage{algorithmic}
\usepackage{xcolor}
\usepackage{colortbl}
\usepackage{tikz,calc}
\usepackage{graphicx}
\usepackage{rotating}
\usepackage{empheq}
\usepackage{blkarray}
\usepackage{eurosym}
\usepackage{amsmath,amsfonts,amssymb,stmaryrd,wasysym}
\usepackage{algorithm,algorithmic}

\usepackage{empheq}
\definecolor{DodgerBlue}{rgb}{0.12, 0.56, 1.0}

\title*{Aitken-Schwarz heterogeneous Domain Decomposition for  EMT-TS Simulation}
\titlerunning{} 
\author{ H. Shourick \inst{2}\inst{1}, D. Tromeur-Dervout \inst{1} and L. Chedot \inst{2}}
\institute{ \inst{1} University of Lyon, UMR5208 U.Lyon1-CNRS,  Institut Camille Jordan, \\ \email{damien.tromeur-dervout@univ-lyon1.fr}\\
 \inst{2} Supergrid-Institute, 14 rue Cyprien,  69200 Villeurbanne. \\ \email{h.shourick,l.chedot@supergrid-institute.com}}
\begin{document}

\maketitle
~\\

\abstract{In this paper, a Schwarz heterogeneous domain decomposition method (DDM) is used to co-simulate an RLC electrical  circuit where a part of the domain is modeled with Electro-Magnetic Transients (EMT) modeling and the other part with dynamic phasor (TS) modeling. Domain partitioning is not based on cutting at transmission lines which introduces a physical delay on the dynamics of the solution, as is usually done, but only on connectivity considerations.
We show the convergence property of the homogeneous  DDM EMT-EMT and TS-TS and of the heterogeneous DDM TS-EMT, with and without overlap and we use the pure linear divergence/convergence  of the method to accelerate it toward the searched solution with the Aitken's acceleration of the convergence  technique. 
}
\vspace*{-0.5cm}
\section{Introduction \label{shourick_contrib_Sec1}}

The introduction of renewable energies into the power grid leads to the use of more components based on power electronics which have to be well dimensioned in order not to be damaged by electrical disturbances. These components imply faster dynamics, for power system safety simulations, which cannot be handled by traditional Transient Simulations (TS) with dynamic phasors. Nevertheless, for large power grids, it can be expected that the need of high level details requiring Electro-Magnetic Transient (EMT) modeling will be localized close to disturbances, as other parts of the network still use TS modeling. 
This paper deals with a proof of concept to develop heterogeneous Schwarz domain decomposition with different modeling (EMT-TS) between the sub-domains. Hybrid (Jacobi type)  EMT-TS co-simulation has to face several locks \cite{shourick_contrib_Interfacing}: EMT and TS do not use the same time step size, the transmission of values is also a problem as the solutions do not have the same representation and are subject to some information loss. Our approach don't use waveform relaxation \cite{shourick_contrib_lelarasmee}, and the domain partitioning is not based on cutting the transmission lines \cite{shourick_contrib_4663663,shourick_contrib_plumwaveform} as we want to be able to define an overlap between the two representations. On the  contrary, we want to use the traditional Schwarz DDM but also where the transmission conditions can lead to divergent DDM. The pure linear convergence/divergence of the linearized problems is then used to accelerate the convergence to the solution by the Aitken's technique.
In Section \ref{shourick_contrib_Sec2}, we describe the EMT and TS modeling and perform homogeneous Schwarz DDM accelerated by the Aitken's acceleration of the convergence technique.  Section \ref{shourick_contrib_Sec3} gives  behavior results obtained for each modeling.
Section \ref{shourick_contrib_Sec4} describes the heterogeneous EMT-TS DDM  and 
gives  first results obtained before concluding in section \ref{shourick_contrib_Sec5}
\vspace{-0.5cm}
\section{EMT and TS modeling \label{shourick_contrib_Sec2}}
Simulation of power grid consists in solving a system of differential algebraic equations (DAE) where the unknowns are currents and voltages. This system is built using the Modified Augmented Nodal Analysis where each component of the grid contributes through  relations   between currents and voltages and the Kirshoff's laws give the algebraical constraints. Let $x$ (respectively $y$) be the differential  (respectively algebraical) unknowns. For the EMT modeling, we have to solve the DAE:
\begin{equation}
F(t,x(t),\dot{x}(t),y(t))=0, \textrm{~with Initial Conditions}  \label{shourick_contrib_DAE}
\end{equation}

The linearized BDF time discretization of \eqref{shourick_contrib_DAE} (Backward Euler here) leads to solve the linear system \eqref{shourick_contrib_LS} to integrate the state space representation of the DAE from time step $t^{n}$ to time step $t^{n+1}$:

\begin{eqnarray}
\underbrace{\left( \begin{array}{cc}
I  -\Delta t A &  B\\
C & D \end{array} \right)}_{\mathbf{H_{\Delta t}}}\left( \begin{array}{c}
x^{n+1} \\ y^{n+1} \end{array}\right) &=&\left( \begin{array}{cc}
I   &  0\\
0 & 0 \end{array} \right) \left( \begin{array}{c}
x^{n} \\ y^{n} \end{array}\right)  \label{shourick_contrib_LS}
\end{eqnarray}
For TS modeling the variables are supposed to oscillate with a specific frequency $\omega_0$ and its selected harmonics taken in a subset $I=\left\{\ldots,-1,0,1, \ldots\right\}$:

\begin{equation}
z(t)= \sum_{k \in I} z_k(t) e^{i k \omega_0 t},\, z=\left\{x,y\right\}.
 \label{shourick_contrib_TS}
\end{equation}
 Introducing \eqref{shourick_contrib_TS} into \eqref{shourick_contrib_DAE} leads after simplification to an another  DAE system that  takes into account the differential property of the dynamic phasor. The resulting DAE system has smoother dynamics. The number of TS variables is then multiplied by the number of harmonics chosen, and the number of equations must be multiplied accordingly. For example,  below is the structure of the matrix $H_{TS}$ by choosing two harmonics $k=a$ and $k=b$ and by solving the imaginary and real part separately and  with $\mathbf{S}$  the matrix  taking into account the differential property of the dynamic phasor modeling.

\begin{eqnarray}
 \mathbf{H}_{TS}\color{black}=
\left(
\begin{array}{c|c}
\begin{array}{c|c}
 \mathbf{H}_{\Delta T}  & -a\,\omega_0 \, \mathbf{S}   \\
\hline
a\,\omega_0 \,\mathbf{S}& \mathbf{H}_{\Delta T} 
\end{array}
&
{\large
0 }\\
& \\
\hline
& \\
0
&
\begin{array}{c|c}
 \mathbf{H}_{\Delta T}  & -c\,\omega_0\,  \mathbf{S} \color{black}  \\
\hline
c\,\omega_0 \, \mathbf{S} & \mathbf{H}_{\Delta T} 
\end{array} 
\end{array}  \right) \nonumber 
\end{eqnarray}
Let  $x_{T}^{n+1}$ (respectively  $x_{E}^{n+1}$) be the algebraic and differential unknowns of TS (respectively EMT) modeling associated to the linear system $H_{TS}x_{T}^{n+1}= b_{T}^n$ (respectively $H_{E}x_{E}^{n+1}= b_{E}^n$)
\vspace{-0.5cm}
\section{EMT and TS Schwarz homogeneous DDM  \label{shourick_contrib_Sec3}}
\begin{figure}[h]
\begin{minipage}{10.5cm}
\begin{minipage}{5.2cm}
 \begin{tikzpicture}[scale=0.6]
 \draw[black!70] (-0.7,3) node[above]{$\Omega$};
 \draw (-4,2) -- (-3,2);
 \draw (-2,2) -- (-1.4,2);
 \draw(0,2) -- (1.5,2);
 \draw (2,2) -- (3,2);
 \draw (3,2) -- (3,-2); 
 \draw (-4,-2) -- (-2.5,-2);
 \draw (-2,-2) -- (-1,-2);
 \draw (0,-2) -- (1,-2);
 \draw (2.4,-2) -- (3,-2);
 \draw (-4,2) -- (-4,-2);

%

  \draw[DodgerBlue] (-4,2) node[above]{\scriptsize $2$};
   \draw[DodgerBlue] (-1.6,2) node[above]{\scriptsize $3$};
   \draw[DodgerBlue] (0.6,2) node[above]{\scriptsize $4$};
   \draw[DodgerBlue] (3,2) node[above]{\scriptsize $5$};
   \draw[DodgerBlue] (-1.6,-2) node[below]{\scriptsize $7$};
   \draw[DodgerBlue] (0.6,-2) node[below]{\scriptsize $6$};
   \draw[DodgerBlue] (-4,-2) node[left]{\scriptsize $1$};
  
 \draw[thick] (1.5,2.5) -- (1.5,1.5);
 \draw[thick] (2,2.5) -- (2,1.5);

  \draw (1.7,2.4) node[above]{\footnotesize  $C_1$};
   
 \draw[thick] (-2,-2.5) -- (-2,-1.5);
 \draw[thick] (-2.5,-2.5) -- (-2.5,-1.5);
 
  \draw (-2.3,-1.6) node[above]{\footnotesize $C_2$};

  \draw[ thick] (-4,-2) -- (-4,-2.5);
 \draw[thick] (-4.5,-2.5) -- (-3.6,-2.5);
  \draw[thick] (-4.5,-2.5) -- (-4.3,-2.7);
 \draw[thick] (-4.2,-2.5) -- (-4,-2.7);
  \draw[thick] (-3.9,-2.5) -- (-3.7,-2.7);
   \draw[thick] (-3.6,-2.5) -- (-3.4,-2.7);
   
   
    \draw[thick] (-1.4,2) -- (-1.2,2.3);
    \draw[thick] (-1.2,2.3) -- (-1,1.7);
     \draw[thick] (-1,1.7) -- (-0.8,2.3);
 \draw[thick] (-0.8,2.3) -- (-0.6,1.7);
  \draw[thick] (-0.6,1.7) -- (-0.4,2.3);
  \draw[thick] (-0.4,2.3) -- (-0.2,1.7);
  \draw[thick] (-0.2,1.7) -- (0,2);
  
   \draw (-0.7,2.4) node[above]{\footnotesize $R_1$};
  
  \draw[thick] (1,-2) -- (1.2,-1.7);
  \draw[thick] (1.2,-1.7) -- (1.4,-2.3);
 \draw[thick] (1.4,-2.3) -- (1.6,-1.7);
 \draw[thick] (1.6,-1.7) -- (1.8,-2.3);
  \draw[thick] (1.8,-2.3) -- (2,-1.7);
  \draw[thick] (2,-1.7) -- (2.2,-2.3);
  \draw[thick] (2.2,-2.3) -- (2.4,-2);
   \draw (1.7,-1.6) node[above]{\footnotesize $R_2$};
  
  \draw[thick] (-4,0) circle(0.5);
  \draw[thick] (-4,0.5)-- (-4,-0.5);
  \draw[->,thick] (-3.45,-0.5)-- (-3.45,0.5);
  \draw (-3.45,0) node[right]{\footnotesize E cos $\omega t = \beta$};

   
   \begin{scope}[shift={(-3.5,2)},rotate=90]
{
 \draw[black!5!yellow!10!red!8!,opacity=0.3] (-0.2,0) rectangle (0.4,-1.5);
 \foreach \r in {0,...,2}
 {
  \draw[thick,scale=1/3,shift={(0,-\r)}]
	(0,0) .. controls ++(2,0) and ++(1,0) ..
	++(0,-1.5) .. controls ++(-1,0) and ++(-0.5,0) ..
	++(0,0.5);
 }
 \draw[thick,scale=1/3,shift={(0,-3)}] (0,0) .. controls ++(2,0) and ++(1,0) .. ++(0,-1.5);
}

\end{scope}
\draw (-2.8,2.4) node[above]{\footnotesize  $L_1$};

   \begin{scope}[shift={(-1.2,-2)},rotate=90]
{
 \draw[black!7!DodgerBlue!11,opacity=0.80] (-0.2,0) rectangle (0.4,-1.5);
 \foreach \r in {0,...,2}
 {
  \draw[thick,scale=1/3,shift={(0,-\r)}]
	(0,0) .. controls ++(2,0) and ++(1,0) ..
	++(0,-1.5) .. controls ++(-1,0) and ++(-0.5,0) ..
	++(0,0.5);
 }
 \draw[thick,scale=1/3,shift={(0,-3)}] (0,0) .. controls ++(2,0) and ++(1,0) .. ++(0,-1.5);
}
\end{scope}
 \draw (-0.5,-1.6) node[above]{\footnotesize $L_2$};

\coordinate (h) at (0.3,2.1);;
\coordinate (i) at (3.5,2.1);
\coordinate (j) at (3.5,-2.2);
\coordinate (k) at (-1.5,-2.2);
\coordinate (hk) at (-0.4,0.15);

\coordinate (m) at (-4.4,2.2);
\coordinate (mn) at (-2,2.9);
\coordinate (n) at (0.2,2.2);
\coordinate (no) at (-0.46,-0.32);
\coordinate (o) at (-1.5,-2.1);
\coordinate (p) at (-4.4,-2.15);
\coordinate (pm) at (-4.85,0);

\draw[dotted,draw=red!60!yellow!35!black,opacity=0.55]  (m) .. controls +(1,0.5) and +(-1,0.08) .. (mn)
               .. controls +(1,0.08) and +(-0.8,0.5) .. (n)
               .. controls +(0.5,-0.8) and +(-0.19,2) .. (no)
               .. controls +(-0.15,-1.5) and +(0.55,0.72) .. (o)
               .. controls +(-1,-0.2) and +(1,-0.2) .. (p)
                .. controls +(-0.5,0.8) and +(0.1,-0.8) .. (pm)
               .. controls +(0.1,0.8) and +(-0.5,-0.8) .. (m); 

  \draw[dotted,draw=black,opacity=0.5](h) .. controls +(0.8,0.6) and +(-0.8,0.6) .. (i)
               .. controls +(0.8,-0.9) and +(0.7,0.9) .. (j)
               .. controls +(-0.9,-0.8) and +(0.9,-0.7) .. (k)
                .. controls +(0.87,0.8) and +(0,-1.3) .. (hk)
                .. controls +(0,1.3) and +(-0.1,-1) .. (h); 
\end{tikzpicture}
\vspace*{-1.cm}
\begin{eqnarray}
v_1&=&0 \\
v_2-v_1-E-Z_s i_{12} &=& 0 \\
v_3-v_2-L_1 \dfrac{d{i}_{23}}{dt} &=& 0 \\
v_4-v_3 -R_1 i_{34} &=& 0 
\end{eqnarray}
\end{minipage}
\hfill
\begin{minipage}{5.4cm}
\begin{eqnarray}
C_1 (\dfrac{d{v}_5}{dt}-\dfrac{d{v}_4}{dt})-i_{45} &=&0\\
v_6-v_5-R_2 i_{56} &=& 0\\
v_7-v_6-L_2 \dfrac{d{i}_{67}}{dt} &=& 0\\
C_2 (\dfrac{d{v}_1}{dt}-\dfrac{d{v}_7}{dt})-i_{71} &=& 0  \\
i_{12}-i_{23}&=& 0  \\
i_{23}-i_{34}&=& 0  \\
i_{34}-i_{45}&=& 0  \\
i_{45}-i_{56}&=& 0  \\
i_{56}-i_{67}&=& 0  \\
i_{67}-i_{71}&=& 0  
\end{eqnarray}
\end{minipage}
\end{minipage}
\caption{Linear RLC circuit and its associated  EMT modeling  DAE system with\\ $x=\left\{ v_1, i_{23}, v4, v5, i_{67}, v_7\right\}$ and $y=\left\{ v_2, i_{12},  v_3, i_{34}, i_{45}, i_{56}, v_6, i_{71}\right\}$. $L1=L2=0.7$, \\ $C1=C2=1.10^{-6}$, $R1=R2=77$, $Zs=1.10^{-6}$, $\omega_0=2\pi\,50$. \label{shourick_contrib_Fig1}}
\end{figure}
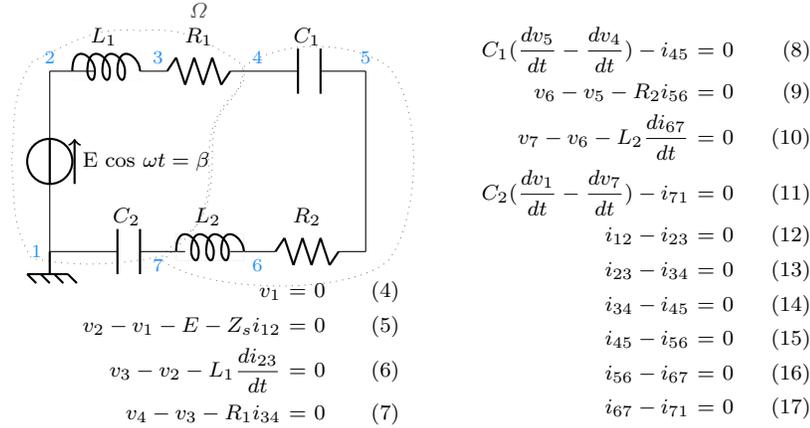

We consider a linear RLC circuit of Figure \ref{shourick_contrib_Fig1} to develop the proof of concept of the  the Schwarz DDM on TS  and EMT models. 

One Restrictive Additive Schwarz (RAS) iterate   to solve  $H x^{m+1,\infty}=b^m \in \mathbb{R}^n$  writes on subdomain $\Omega_i$:  $x_i^{m+1,k+1} = A_i^{-1} \left(b_i^m - E_i x_{i,e}^{m+1,k} \right)$, with $R_i \in \mathbb{R}^{n_i \times n}$ the operator that restricts the global vector  to the  subdomain $\Omega_i$, including the overlap,  $\widetilde{R}_i\in \mathbb{R}^{n_i \times n}$ the operator that restricts the global vector to the subdomain  $\Omega_i$, with setting to $0$ the components of the vector that correspond to the overlap. $W_i$ is the global index set  of the unknowns belonging to the subdomain $\Omega_i$.  $A_i$ is the part of the operator $H$ associated to the subdomain $\Omega_i$: $A_i = R_i H R_i^T$.  $x_i^{m+1} = R_i x^{m+1}$ and $b_i^m = R_i b^m$ are the restriction to the subdomain $\Omega_i$ of the solution and the right hand side respectively. $x_{i,e}^{m+1}$ represents the external data dependencies  of the subdomain $\Omega_i$ : $x_{i,e}^{m+1}$  is composed of the  $x_j^{m+1}$ such that $H_{kj} \neq 0$ with $k \in W_i$ and $j \notin W_i$.  $R_{i,e}$ is the restriction operator such that  $R_{i,e} x^{m+1} = x_{i,e}^{m+1}$. $E_i$  is the part of the matrix $H$ that represents  the effect of the unknowns external to the subdomain $\Omega_i$ on the unknowns belonging to the subdomain $\Omega_i$ : $E_i  = R_{i,e} H R_{i,e}^T$.

\begin{figure}[h]
\centering
\begin{minipage}{11cm}

\begin{minipage}{5.4cm}
\begin{tikzpicture}[scale=0.5]

 \draw[black!70] (-0.7,3) node[above]{$\Omega = \Omega_1 \cup \Omega_2$};
 \draw[black!60!DodgerBlue!70] (2.55,2.7) node[above]{$\Omega_2$};
 \draw[black!70!yellow!65!red!65!] (-3.65,2.7) node[above]{$\Omega_1$};
 \draw (-4,2) -- (-3,2);
 \draw (-2,2) -- (-1.4,2);
 \draw(0,2) -- (1.5,2);
 \draw (2,2) -- (3,2);
 \draw (3,2) -- (3,-2); 
 \draw (-4,-2) -- (-2.5,-2);
 \draw (-2,-2) -- (-1,-2);
 \draw (0,-2) -- (1,-2);
 \draw (2.4,-2) -- (3,-2);
 
  \draw (-4,2) -- (-4,-2);

   \fill[DodgerBlue] (-4,2) circle(0.05);
   \fill[DodgerBlue] (-1.6,2) circle(0.05);
   \fill[DodgerBlue] (0.6,2) circle(0.05);
   \fill[DodgerBlue] (3,2) circle(0.05);
   \fill[DodgerBlue] (-1.6,-2) circle(0.05);
   \fill[DodgerBlue] (0.6,-2) circle(0.05);
   \fill[DodgerBlue] (-4,-2) circle(0.05);

  \draw[DodgerBlue] (-4,2) node[above]{\scriptsize $2$};
   \draw[DodgerBlue] (-1.6,2) node[above]{\scriptsize $3$};
   \draw[DodgerBlue] (0.6,2) node[above]{\scriptsize $4$};
   \draw[DodgerBlue] (3,2) node[above]{\scriptsize $5$};
   \draw[DodgerBlue] (-1.6,-2) node[below]{\scriptsize $7$};
   \draw[DodgerBlue] (0.6,-2) node[below]{\scriptsize $6$};
   \draw[DodgerBlue] (-4,-2) node[left]{\scriptsize $1$};
  
 \draw[thick] (1.5,2.5) -- (1.5,1.5);
 \draw[thick] (2,2.5) -- (2,1.5);

  \draw (1.7,2.4) node[above]{\footnotesize  $C_1$};
   
 \draw[thick] (-2,-2.5) -- (-2,-1.5);
 \draw[thick] (-2.5,-2.5) -- (-2.5,-1.5);
 
  \draw (-2.3,-1.6) node[above]{\footnotesize $C_2$};

  \draw[ thick] (-4,-2) -- (-4,-2.5);
 \draw[thick] (-4.5,-2.5) -- (-3.6,-2.5);
  \draw[thick] (-4.5,-2.5) -- (-4.3,-2.7);
 \draw[thick] (-4.2,-2.5) -- (-4,-2.7);
  \draw[thick] (-3.9,-2.5) -- (-3.7,-2.7);
   \draw[thick] (-3.6,-2.5) -- (-3.4,-2.7);
   
   
    \draw[thick] (-1.4,2) -- (-1.2,2.3);
    \draw[thick] (-1.2,2.3) -- (-1,1.7);
     \draw[thick] (-1,1.7) -- (-0.8,2.3);
 \draw[thick] (-0.8,2.3) -- (-0.6,1.7);
  \draw[thick] (-0.6,1.7) -- (-0.4,2.3);
  \draw[thick] (-0.4,2.3) -- (-0.2,1.7);
  \draw[thick] (-0.2,1.7) -- (0,2);
  
   \draw (-0.7,2.4) node[above]{\footnotesize $R_1$};
  
  \draw[thick] (1,-2) -- (1.2,-1.7);
  \draw[thick] (1.2,-1.7) -- (1.4,-2.3);
 \draw[thick] (1.4,-2.3) -- (1.6,-1.7);
 \draw[thick] (1.6,-1.7) -- (1.8,-2.3);
  \draw[thick] (1.8,-2.3) -- (2,-1.7);
  \draw[thick] (2,-1.7) -- (2.2,-2.3);
  \draw[thick] (2.2,-2.3) -- (2.4,-2);
   \draw (1.7,-1.6) node[above]{\footnotesize $R_2$};
  
  \draw[thick] (-4,0) circle(0.5);
  \draw[thick] (-4,0.5)-- (-4,-0.5);
  \draw[->,thick] (-3.45,-0.5)-- (-3.45,0.5);
  \draw (-3.45,0) node[right]{\footnotesize E cos $\omega t = \beta$};

   
   \begin{scope}[shift={(-3.5,2)},rotate=90]
{
 \draw[black!5!yellow!10!red!8!,opacity=0.3,fill=black!5!yellow!10!red!8!,opacity=0.8] (-0.2,0) rectangle (0.4,-1.5);
 \foreach \r in {0,...,2}
 {
  \draw[thick,scale=1/3,shift={(0,-\r)}]
	(0,0) .. controls ++(2,0) and ++(1,0) ..
	++(0,-1.5) .. controls ++(-1,0) and ++(-0.5,0) ..
	++(0,0.5);
 }
 \draw[thick,scale=1/3,shift={(0,-3)}] (0,0) .. controls ++(2,0) and ++(1,0) .. ++(0,-1.5);
}

\end{scope}
\draw (-2.8,2.4) node[above]{\footnotesize  $L_1$};

   \begin{scope}[shift={(-1.2,-2)},rotate=90]
{
 \draw[black!7!DodgerBlue!11,opacity=0.80,fill=black!7!DodgerBlue!11,opacity=0.80] (-0.2,0) rectangle (0.4,-1.5);
 \foreach \r in {0,...,2}
 {
  \draw[thick,scale=1/3,shift={(0,-\r)}]
	(0,0) .. controls ++(2,0) and ++(1,0) ..
	++(0,-1.5) .. controls ++(-1,0) and ++(-0.5,0) ..
	++(0,0.5);
 }
 \draw[thick,scale=1/3,shift={(0,-3)}] (0,0) .. controls ++(2,0) and ++(1,0) .. ++(0,-1.5);
}
\end{scope}
 \draw (-0.5,-1.6) node[above]{\footnotesize $L_2$};

\draw[red,thick,->] (0.1,2.2) -- (0.5,2.2);
\draw[blue,thick,<-] (-1.4,-2.1) -- (-1,-2.1);
\draw[red](-0.1,2.1) node[above]{\footnotesize ${\bf i_{34},v_3}$};
\draw[blue](-0.8,-1.9) node[below]{\footnotesize ${\bf i_{67},v_6}$};

\coordinate (h) at (0.3,2.1);;
\coordinate (i) at (3.5,2.1);
\coordinate (j) at (3.5,-2.2);
\coordinate (k) at (-1.5,-2.2);
\coordinate (hk) at (-0.4,0.15);

\coordinate (m) at (-4.4,2.2);
\coordinate (mn) at (-2,2.9);
\coordinate (n) at (0.2,2.2);
\coordinate (no) at (-0.46,-0.32);
\coordinate (o) at (-1.5,-2.1);
\coordinate (p) at (-4.4,-2.15);
\coordinate (pm) at (-4.85,0);

\draw[dotted,draw=red!60!yellow!35!black,fill=black!5!yellow!10!red!15!,opacity=0.55]  (m) .. controls +(1,0.5) and +(-1,0.08) .. (mn)
               .. controls +(1,0.08) and +(-0.8,0.5) .. (n)
               .. controls +(0.5,-0.8) and +(-0.19,2) .. (no)
               .. controls +(-0.15,-1.5) and +(0.55,0.72) .. (o)
               .. controls +(-1,-0.2) and +(1,-0.2) .. (p)
                .. controls +(-0.5,0.8) and +(0.1,-0.8) .. (pm)
               .. controls +(0.1,0.8) and +(-0.5,-0.8) .. (m); 

  \draw[dotted,draw=black,fill=black!5!DodgerBlue!11,opacity=0.5](h) .. controls +(0.8,0.6) and +(-0.8,0.6) .. (i)
               .. controls +(0.8,-0.9) and +(0.7,0.9) .. (j)
               .. controls +(-0.9,-0.8) and +(0.9,-0.7) .. (k)
                .. controls +(0.87,0.8) and +(0,-1.3) .. (hk)
                .. controls +(0,1.3) and +(-0.1,-1) .. (h); 

\end{tikzpicture}
\end{minipage}
\hfill
\begin{minipage}{5.4cm}
\includegraphics[scale=0.20]{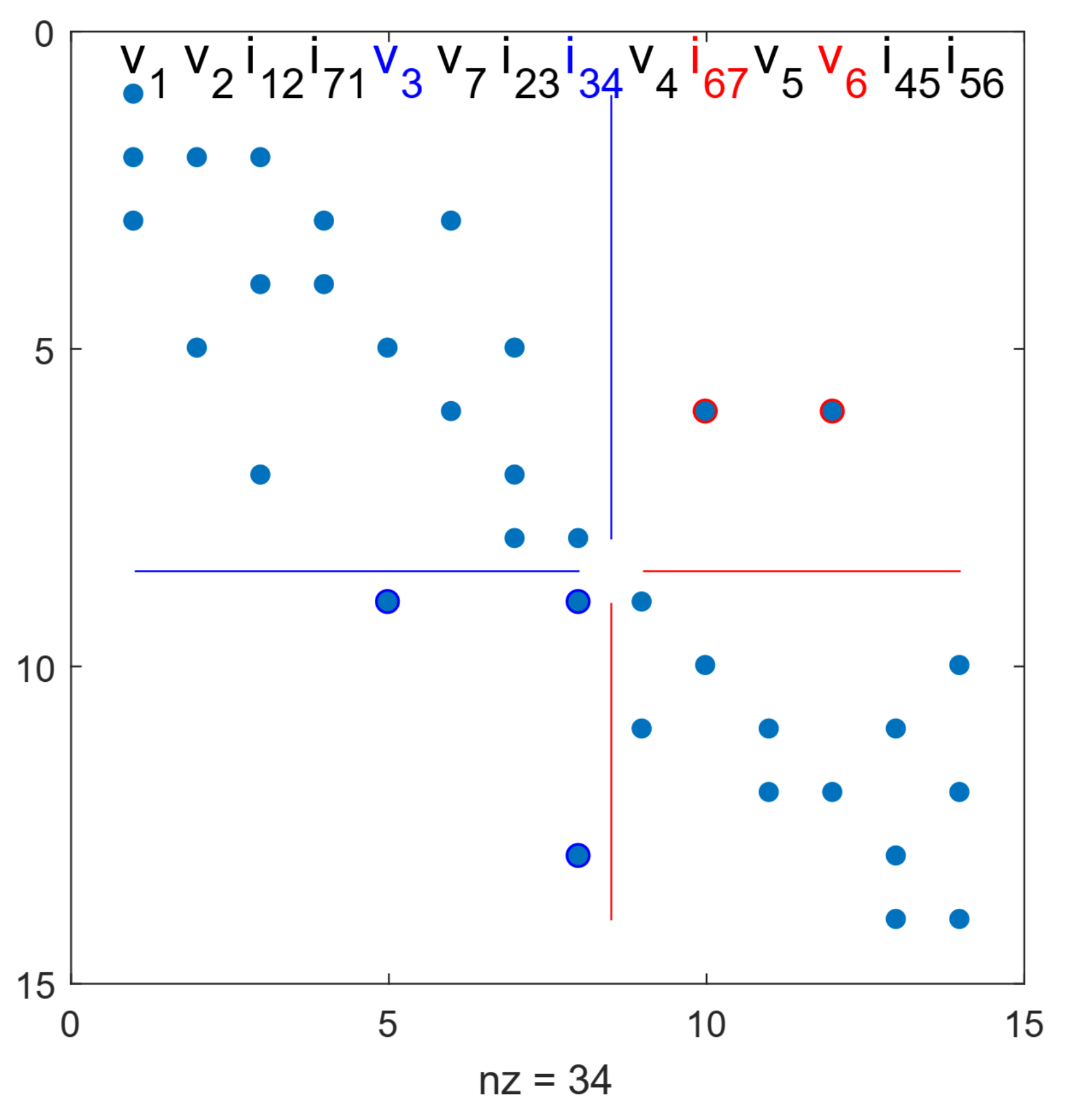}
\end{minipage}
\end{minipage}
\begin{minipage}{11cm}

\begin{minipage}{5.4cm}

 \begin{tikzpicture}[scale=0.5]

 \draw[black!70] (-0.7,3) node[above]{$\Omega = \Omega_1 \cup \Omega_2$};
 \draw[black!60!DodgerBlue!70] (2.55,2.7) node[above]{$\Omega_2$};
 \draw[black!70!yellow!65!red!65!] (-3.65,2.7) node[above]{$\Omega_1$};
 \draw (-4,2) -- (-3,2);
 \draw (-2,2) -- (-1.4,2);
 \draw(0,2) -- (1.5,2);
 \draw (2,2) -- (3,2);
 \draw (3,2) -- (3,-2); 
 \draw (-4,-2) -- (-2.5,-2);
 \draw (-2,-2) -- (-1,-2);
 \draw (0,-2) -- (1,-2);
 \draw (2.4,-2) -- (3,-2);
 
  \draw (-4,2) -- (-4,-2);

   \fill[DodgerBlue] (-4,2) circle(0.05);
   \fill[DodgerBlue] (-1.6,2) circle(0.05);
   \fill[DodgerBlue] (0.6,2) circle(0.05);
   \fill[DodgerBlue] (3,2) circle(0.05);
   \fill[DodgerBlue] (-1.6,-2) circle(0.05);
   \fill[DodgerBlue] (0.6,-2) circle(0.05);
   \fill[DodgerBlue] (-4,-2) circle(0.05);

  \draw[DodgerBlue] (-4,2) node[above]{\scriptsize $2$};
   \draw[DodgerBlue] (-1.6,2) node[above]{\scriptsize $3$};
   \draw[DodgerBlue] (0.6,2) node[above]{\scriptsize $4$};
   \draw[DodgerBlue] (3,2) node[above]{\scriptsize $5$};
   \draw[DodgerBlue] (-1.6,-2) node[below]{\scriptsize $7$};
   \draw[DodgerBlue] (0.6,-2) node[below]{\scriptsize $6$};
   \draw[DodgerBlue] (-4,-2) node[left]{\scriptsize $1$};
  
 \draw[thick] (1.5,2.5) -- (1.5,1.5);
 \draw[thick] (2,2.5) -- (2,1.5);

  \draw (1.7,2.4) node[above]{\footnotesize  $C_1$};
   
 \draw[thick] (-2,-2.5) -- (-2,-1.5);
 \draw[thick] (-2.5,-2.5) -- (-2.5,-1.5);
 
  \draw (-2.3,-1.6) node[above]{\footnotesize $C_2$};

  \draw[ thick] (-4,-2) -- (-4,-2.5);
 \draw[thick] (-4.5,-2.5) -- (-3.6,-2.5);
  \draw[thick] (-4.5,-2.5) -- (-4.3,-2.7);
 \draw[thick] (-4.2,-2.5) -- (-4,-2.7);
  \draw[thick] (-3.9,-2.5) -- (-3.7,-2.7);
   \draw[thick] (-3.6,-2.5) -- (-3.4,-2.7);
   
   
    \draw[thick] (-1.4,2) -- (-1.2,2.3);
    \draw[thick] (-1.2,2.3) -- (-1,1.7);
     \draw[thick] (-1,1.7) -- (-0.8,2.3);
 \draw[thick] (-0.8,2.3) -- (-0.6,1.7);
  \draw[thick] (-0.6,1.7) -- (-0.4,2.3);
  \draw[thick] (-0.4,2.3) -- (-0.2,1.7);
  \draw[thick] (-0.2,1.7) -- (0,2);
  
   \draw (-0.7,2.4) node[above]{\footnotesize $R_1$};
  
  \draw[thick] (1,-2) -- (1.2,-1.7);
  \draw[thick] (1.2,-1.7) -- (1.4,-2.3);
 \draw[thick] (1.4,-2.3) -- (1.6,-1.7);
 \draw[thick] (1.6,-1.7) -- (1.8,-2.3);
  \draw[thick] (1.8,-2.3) -- (2,-1.7);
  \draw[thick] (2,-1.7) -- (2.2,-2.3);
  \draw[thick] (2.2,-2.3) -- (2.4,-2);
   \draw (1.7,-1.6) node[above]{\footnotesize $R_2$};
  
  \draw[thick] (-4,0) circle(0.5);
  \draw[thick] (-4,0.5)-- (-4,-0.5);
  \draw[->,thick] (-3.45,-0.5)-- (-3.45,0.5);
  \draw (-3.45,0) node[right]{\footnotesize E cos $\omega t = \beta$};

   
   \begin{scope}[shift={(-3.5,2)},rotate=90]
{
 \draw[black!5!yellow!10!red!8!,opacity=0.3,fill=black!5!yellow!10!red!8!,opacity=0.8] (-0.2,0) rectangle (0.4,-1.5);
 \foreach \r in {0,...,2}
 {
  \draw[thick,scale=1/3,shift={(0,-\r)}]
	(0,0) .. controls ++(2,0) and ++(1,0) ..
	++(0,-1.5) .. controls ++(-1,0) and ++(-0.5,0) ..
	++(0,0.5);
 }
 \draw[thick,scale=1/3,shift={(0,-3)}] (0,0) .. controls ++(2,0) and ++(1,0) .. ++(0,-1.5);
}

\end{scope}
\draw (-2.8,2.4) node[above]{\footnotesize  $L_1$};

   \begin{scope}[shift={(-1.2,-2)},rotate=90]
{
 \draw[black!7!DodgerBlue!11,opacity=0.80,fill=black!7!DodgerBlue!11,opacity=0.80] (-0.2,0) rectangle (0.4,-1.5);
 \foreach \r in {0,...,2}
 {
  \draw[thick,scale=1/3,shift={(0,-\r)}]
	(0,0) .. controls ++(2,0) and ++(1,0) ..
	++(0,-1.5) .. controls ++(-1,0) and ++(-0.5,0) ..
	++(0,0.5);
 }
 \draw[thick,scale=1/3,shift={(0,-3)}] (0,0) .. controls ++(2,0) and ++(1,0) .. ++(0,-1.5);
}
\end{scope}
 \draw (-0.5,-1.6) node[above]{\footnotesize $L_2$};

\draw[red,thick,->] (-4.5,2.6) -- (-4,2.6);
\draw[blue,thick,<-] (0.9,-1.3) -- (1.3,-1.3);
\draw[red](-4.4,2.45) node[above]{\footnotesize ${\bf i_{12},v_2}$};
\draw[blue](1.3,-1.1) node[below]{\footnotesize ${\bf i_{56},v_6}$};

\coordinate (h) at (-3.8,2.5);
\coordinate (i) at (3.5,2.5);
\coordinate (j) at (3.5,-2.3);
\coordinate (k) at (-1.9,-2.3);

\coordinate (m) at (-4.4,2.2);
\coordinate (mn) at (-1.5,2.8);
\coordinate (n) at (0.63,2.2);
\coordinate (no) at (1,0);
\coordinate (o) at (0.32,-2.2);
\coordinate (op) at (-2,-2.9);
\coordinate (p) at (-4.4,-2.2);
\coordinate (pm) at (-4.85,0);

\draw[dotted,draw=red!60!yellow!35!black,fill=black!5!yellow!10!red!15!,opacity=0.50]  (m) .. controls +(1,0.5) and +(-1,0.08) .. (mn)
               .. controls +(1,0.08) and +(-0.8,0.5) .. (n)
               .. controls +(0.5,-0.8) and +(-0.17,2) .. (no)
               .. controls +(-0.17,-2) and +(0.5,0.8) .. (o)
               .. controls +(-0.7,-0.5) and +(1,0.1) .. (op)
               .. controls +(-1,0.1) and +(0.7,-0.5) .. (p)
                .. controls +(-0.5,0.8) and +(0.1,-0.8) .. (pm)
               .. controls +(0.1,0.8) and +(-0.5,-0.8) .. (m); 

  \draw[dotted,draw=black,fill=black!5!DodgerBlue!11,opacity=0.5](h) .. controls +(1,0.8) and +(-0.9,0.8) .. (i)
               .. controls +(0.8,-0.9) and +(0.7,0.9) .. (j)
               .. controls +(-0.9,-0.8) and +(0.9,-0.7) .. (k)
               .. controls +(-1,1) and +(0.9,-1) .. (h); 

\end{tikzpicture}
\end{minipage}
\hfill
\begin{minipage}{5.4cm}\includegraphics[scale=0.21]{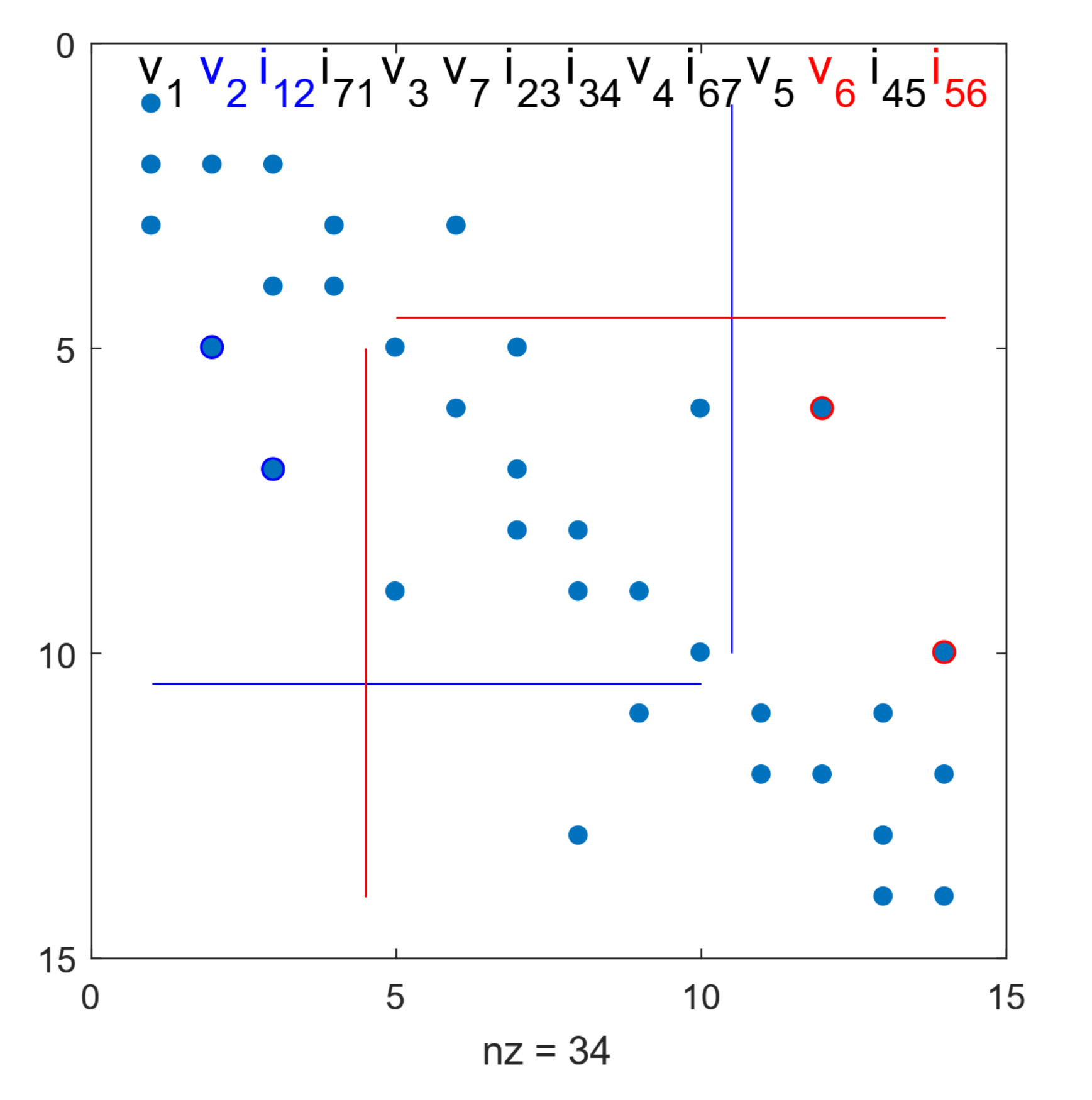}
\end{minipage}
\end{minipage}

\caption{Graph partitionning of the RLC circuit in two subdomains and the associated matrix partioning without overlap (top) and with overlap of 1 (bottom). \label{shourick_contrib_Fig2}}
\end{figure}
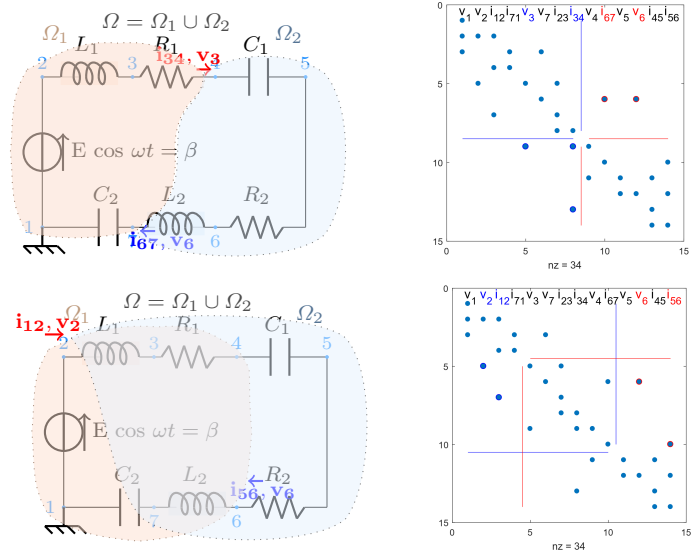

The small linear system  associated with the RLC circuit is partitioned into two  subdomains using graph partitioning  without overlaping (Figure \ref{shourick_contrib_Fig2} top) and with an overlap of 1 (Figure \ref{shourick_contrib_Fig2} bottom). Each subdomain needs two values from the other to solve its equations.

The RAS applied to each time  step has a pure linear convergence i.e. the error operator $P$ does not depend on the RAS iterate.\begin{eqnarray}
x^{m+1,p+1}-x^{m+1,\infty} &=& P (x^{m+1,p}-x^{m+1,\infty}) 
\end{eqnarray}
 Thus it can be accelerated if it does not stagnate to obtain the searched solution regardless of its convergence or divergence \cite{shourick_contrib_dtd}.
\begin{eqnarray}
x^{m+1,\color{black}\infty\color{black}} &=& (I_d-P)^{-1} (x^{m+1,\color{black}1\color{black}}-P x^{m+1,\color{black}0\color{black}})\label{shourick_contrib_EqAitken}
\end{eqnarray}
$P$ can be compute numerically from the values of the iterated transmission conditions. For this small problem it can be directly computed working on the matrix partitioning.
\begin{eqnarray}
P&=&-[(\tilde{R}_1)^tA_{1}^{-1} E_{1,e} R_{1,e}+(\tilde{R}_2)^tA_{2}^{-1} E_{2,e} R_{2,e}] 
\end{eqnarray}
\begin{table}[t]
\centering
\begin{tabular}{|l|c|c|c|c|}
\hline
$\lambda(P)$ & \begin{tabular}{c}without \\ overlap \end{tabular}& \begin{tabular}{c} with \\ overlap \end{tabular}& Schwarz & $\Delta t$\\
\hline
 EMT & $\pm$ 6.0638i &  $ \pm$ 6.0638i  & RAS  & $2.10^{-4}$\\
              \hline
TS {\tiny k=1}  &  -36.6318 $\pm $4.4466i & -36.6318 $\pm $4.4466i  & RMS& $2.10^{-4}$\\
              \hline
TS {\tiny k=0}  &  -36.77 $\pm $0i & -36.77 $\pm $0i $ \pm$ 0i & RMS & $2.10^{-4}$\\
              \hline
           
TS {\tiny k=1}  &  -1.28888$\pm $0.188i & -1.28888$\pm $0.188i & RMS& $2.10^{-3}$\\
              \hline
TS {\tiny k=0}  &  -1.427 $\pm $0i & -1.427 $\pm $0i  & RMS & $2.10^{-3}$\\
              \hline
\end{tabular}
\caption{Larger eigenvalue for  $P$ error operator for RAS and  EMT modeling ($\Delta t=2.10^{-4}$), and for RMS and TS $k=0,1$ ($\Delta T=2.10^{-4}$, $\Delta T=2.10^{-3}$) modeling. \label{shourick_contrib_Tab1}}
\end{table}\\
Table \ref{shourick_contrib_Tab1} gives the larger eigenvalue in modulus for the $P$ RAS error operator for the EMT modeling and for the $P$ RMS(Restricted Multiplicative Schwarz) error operator  for the TS modeling main harmonic $k=1$ applied to the RLC circuit. In both cases EMT and TS modeling the eigenvalue  modulus is greater than one, so the method diverges. We can observe that the overlap does not impact the divergence of the method. The time step increasing from $\Delta t=2.10^{-4}$ to $\Delta T=2.10^{-3}$ has a beneficial effect on the TS-TS DDM divergence. Nevertheless, the divergence is purely linear and the Aitken's acceleration \eqref{shourick_contrib_EqAitken} can be performed after the first iterate.

\begin{figure}[H]
\begin{minipage}{12cm}
\begin{minipage}{5.8cm}
\includegraphics[scale=0.48]{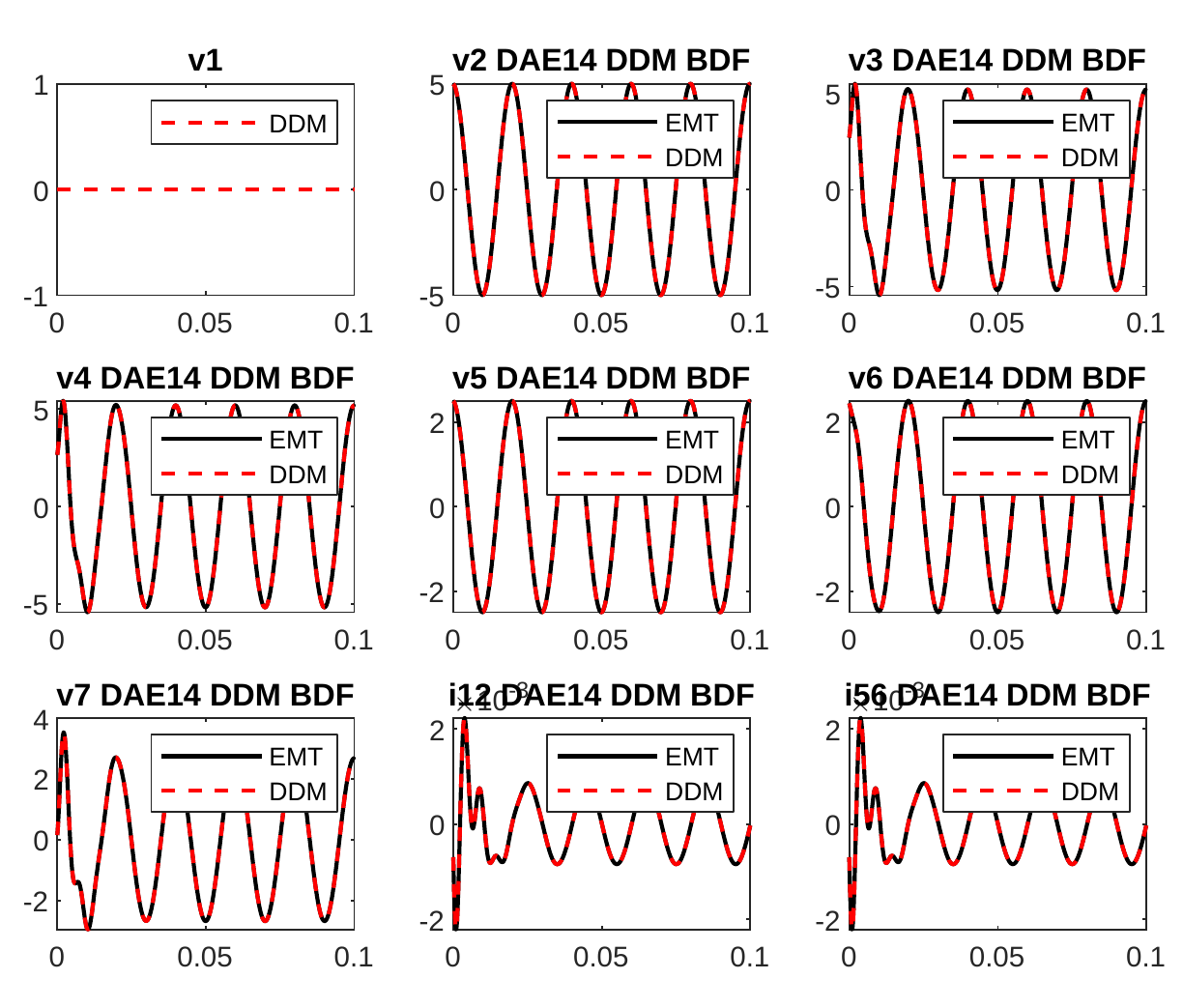}
\end{minipage}
\hfill
\begin{minipage}{5.8cm}
\includegraphics[scale=0.48]{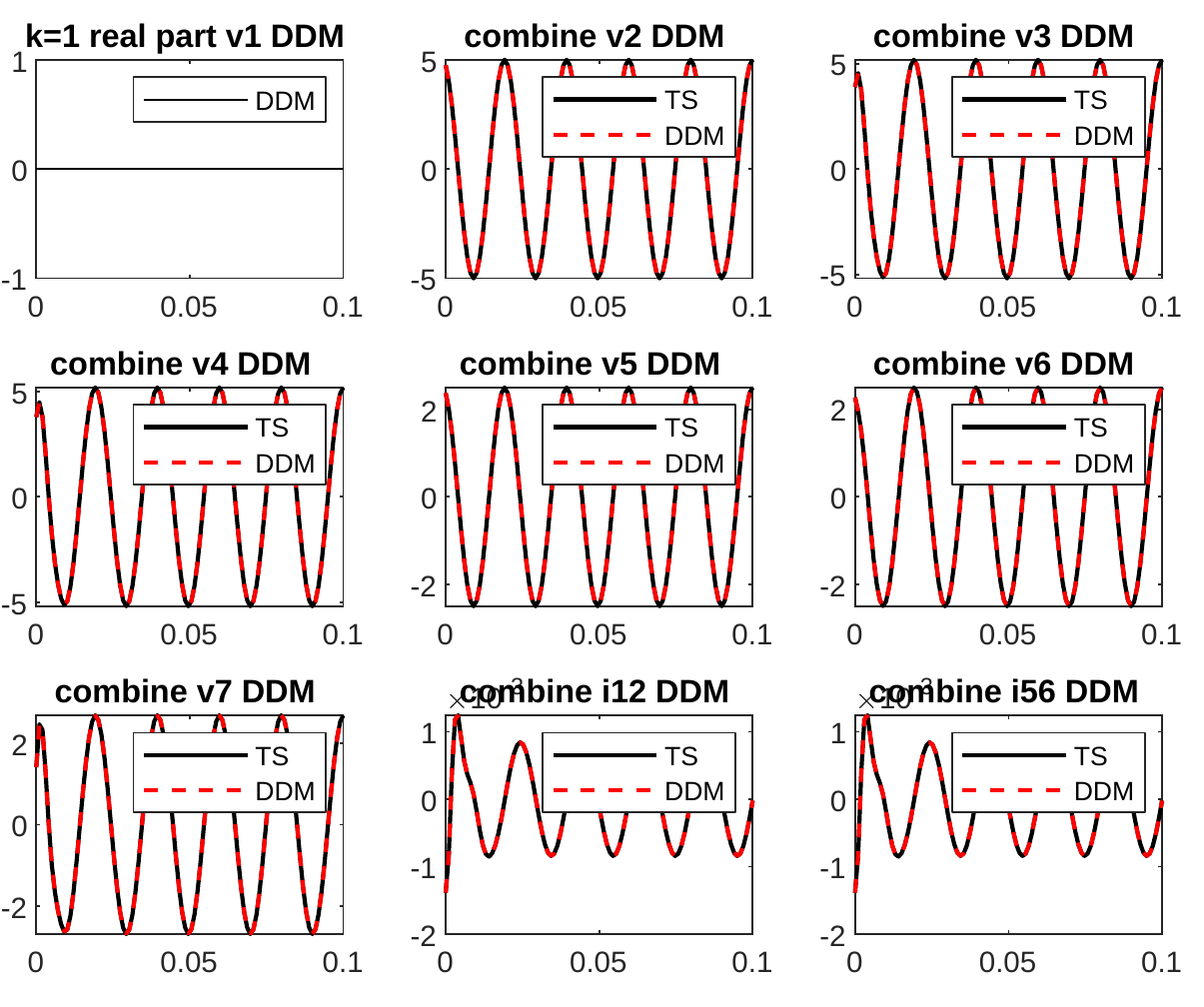}
\end{minipage}
\end{minipage}
\caption{Homogeneous DDM results  comparison with DAE monodomain: (Left) RAS  for EMT modeling with $\Delta t_E=1.10^{-4}$  and (right) RMS for  TS modeling with $\Delta t_T=2.10^{-3}$. \label{shourick_contrib_Fig3}}
\end{figure}

\section{Heterogeneous DDM EMT-TS  \label{shourick_contrib_Sec4}}
Our goal is to simulate, using heterogeneous RAS DDM, the electrical network  with one part with  a TS modeling which can use large time steps $\Delta T$ and the other part with the EMT modeling which requires smaller time steps $\Delta t$ as the high oscillations remain.

These two representations  TS and EMT of the solution imply having some operators $E_{emt}^{TS}$ (respectively $E_{TS}^{emt}$) to transfer  the solution from the subdomain EMT (respectively TS) to the other TS  (respectively EMT). The $E^{emt}_{TS}$ operator needs to compute the fundamental harmonic and other harmonics chosen of the solution from the history of the EMT solution. The history time length is one period. This is performed by the FFT of the solution over the time period and keeping the mode corresponding to the chosen harmonics. 

The $E_{emt}^{TS}$ operator is more simple as it consists in recombining the TS modes of the solution with the appropriate Fourier basis modes.

Let us consider a linear electrical network with the TS modeling. The time discretisation of the DAE to integrate from $T^N$ to $T^{N+1}$, assuming that $\Delta T = m \Delta t$ can be witten as:
\begin{eqnarray}
\underbrace{\left( \begin{array}{cc}
I  -\Delta T A_{TS} &  B_{TS}\\
C_{TS} & D_{TS} \end{array} \right)}_{H_{TS}} \underbrace{\left( \begin{array}{c}
x^{N+1}_{TS} \\ y^{N+1}_{TS} \end{array}\right)}_{w^{N+1}_{TS}} &=& \underbrace{\left( \begin{array}{cc}
I   &  0\\
0 & 0 \end{array} \right)}_{\Theta_{emt}}\left( \begin{array}{c}
x^{N}_{TS} \\ y^{N}_{TS} \end{array}\right)\nonumber\\
&& + \underbrace{\left( \begin{array}{cc}
E^{A}_{TS} &  E^{B}_{TS}\\
E^{C}_{TS} & E^{D}_{TS}
\end{array}\right)}_{E^{emt}_{TS}} \left(\begin{array}{c}
x^{m}_{emt} \\ y^{m}_{emt} \end{array}\right)  
\end{eqnarray}
Similarly one time step for the EMT side to integrate from $t^n$ to $t^{n+1}$ can be witten as:
\begin{eqnarray}
\underbrace{\left( \begin{array}{cc}
I  -\Delta t A_{emt} &  B_{emt}\\
C_{emt} & D_{emt} \end{array} \right)}_{H_{emt}} \underbrace{\left( \begin{array}{c}
x^{n+1}_{emt} \\ y^{n+1}_{emt} \end{array}\right)}_{w^{n+1}} &=&\underbrace{\left( \begin{array}{cc}
I   &  0\\
0 & 0 \end{array} \right)}_{\Theta_{emt}} \left( \begin{array}{c}
x^{n}_{emt} \\ y^{n}_{emt} \end{array}\right)\nonumber\\ 
&& + \underbrace{\left( \begin{array}{cc}
E^{A}_{emt} &  E^{B}_{emt}\\
E^{C}_{emt} & E^{D}_{emt}
\end{array}\right)}_{E_{emt}^{TS}} \underbrace{\left(\begin{array}{c}
x^{N+1}_{TS}(t^{n+1)} \\ y^{N+1}_{TS} (t^{n+1}) \end{array}\right)}_{W^{N+1}(t^{n+1})}
\end{eqnarray}

The $m$ time steps can be gathered in one larger system considering $t^0=T^N$:

\begin{eqnarray}
\underbrace{\left( \begin{array}{cccccc} 
I & & & & &   \\
-\Theta_{emt} & H_{emt} &&&&\\
&-\Theta_{emt} & H_{emt} &&& \\
& & \ddots & \ddots  && \\
& & & -\Theta_{emt} & H_{emt}& \\
&&&& -\Theta_{emt} & H_{emt}  
\end{array}\right)}_{\mathbb{H}_{emt}}\underbrace{\left( \begin{array}{c}
w^0_{emt} \\ w^1_{emt} \\ w^2_{emt} \\ \vdots \\ w^{m-1}_{emt} \\ w^{m}_{emt} 
 \end{array}\right)}_{\mathbb{W}_{emt}} &=& \nonumber \\
 \underbrace{\left( \begin{array}{cccccc} 
I & & & & &   \\
 & E_{emt}^{TS} &&&&\\
& & E_{emt}^{TS} &&& \\
& & \ddots & \ddots  && \\
& & &  & E_{emt}^{TS}& \\
&&&& & E_{emt}^{TS}  
\end{array}\right)}_{\mathbb{E}^{TS}_{emt}}\underbrace{\left( \begin{array}{c}
(x^0,y^0)^t \\ W^{N+1}(t^1) \\ W^{N+1}(t^2) \\ \vdots \\ W^{N+1}(t^{p-1}) \\ W^{N+1}(t^{p}) 
 \end{array}\right)}_{\mathbb{W}^{N+1}_{TS}} 
\end{eqnarray}
This system needs the values that the TS solution connected to the EMT part taken on the small time steps.

The two domains are connected via the connected or flowing variables. Since these variables should be the solution at time $T^{N+1}$, we need the  Schwarz iterative algorithm to obtain the exact values. We then iterate the iteration $p+1$ by taking the connected values, at the iteration $p$, from the other subdomain. We can used the  multiplicative form or the additive form as follows:

\begin{eqnarray}
\left\{ \begin{array}{lcl}
H_{TS} \, w^{N+1,{\bf p+1}}_{TS} &=& \Theta_{TS} \, w^N_{TS} + E_{TS}^{emt}\, w_{emt}^{m,{\bf p}}\\
\mathbb{H}_{emt} \, \mathbb{W}_{emt}^{N+1,\bf p+1} &=& \mathbb{E}_{emt}^{TS} \mathbb{W}^{N+1,{\bf p}}_{TS} \end{array} \right.
\end{eqnarray}

\begin{figure}[t]
\centering
\begin{minipage}{10cm}
\begin{minipage}{3.8cm}
\includegraphics[scale=0.48]{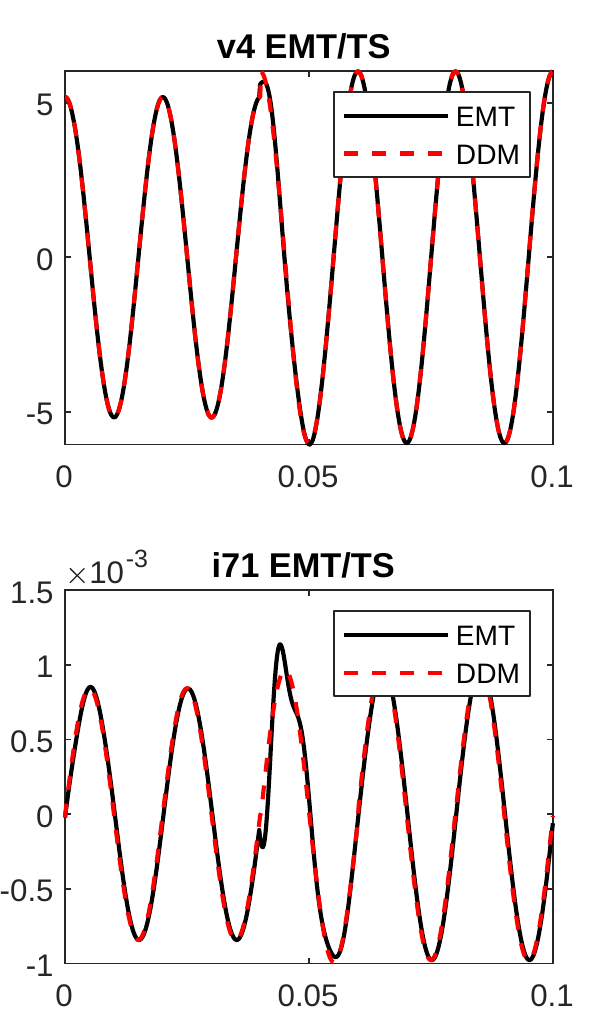}
\end{minipage}
\hfill
\begin{minipage}{5.8cm}
\includegraphics[scale=0.42]{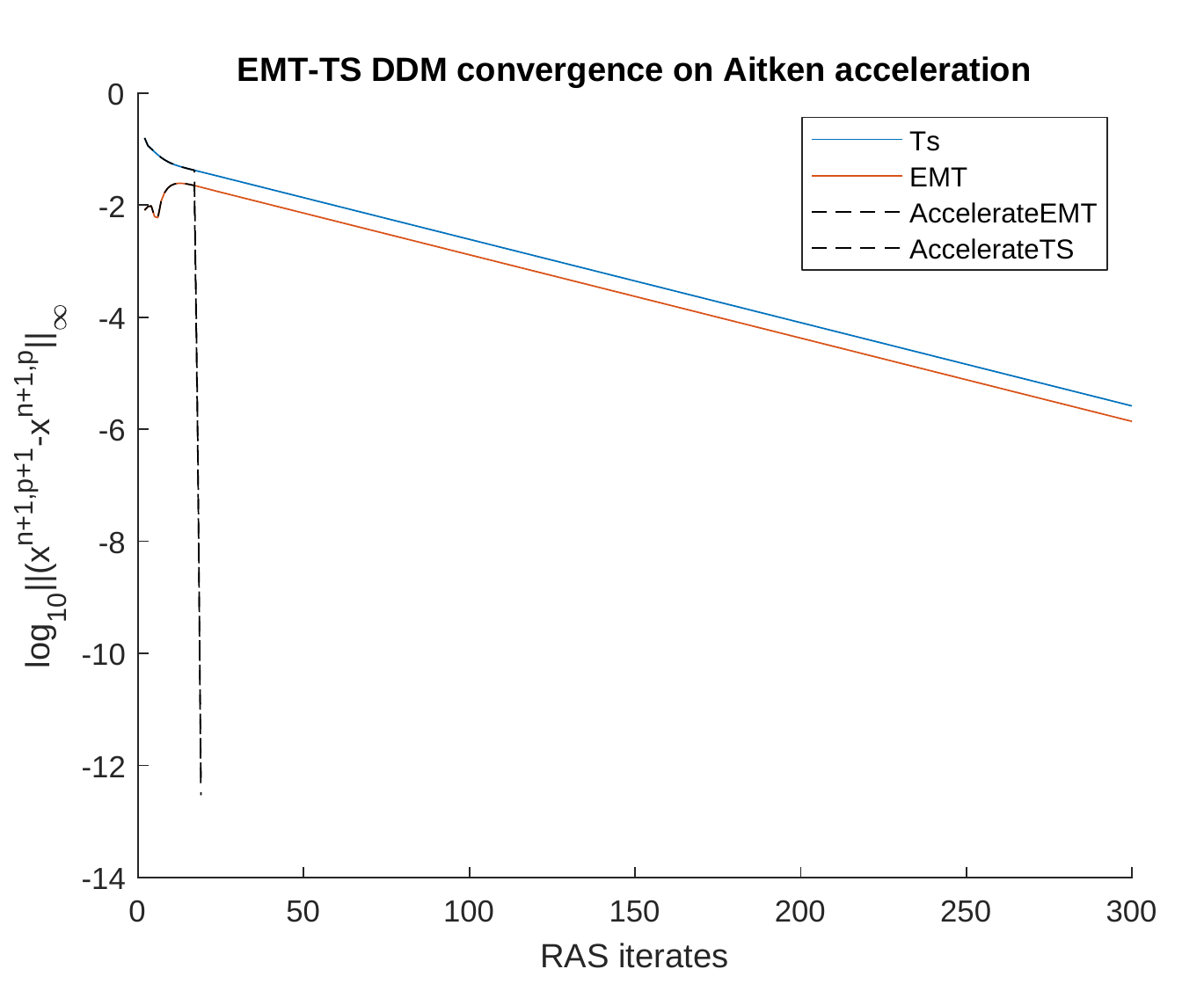}
\end{minipage}
\end{minipage}
\caption{Heterogeneous EMT ($\Delta t=2.10^{-4}$)-TS($\Delta T=2.10^{-2}$) DDM results  comparison with DAE monodomain (Left) and RAS convergence error for each subdomain at $t=0.02$ and its Aitken's acceleration with $P$ computed numerically from 9 iterates (right). \label{shourick_contrib_Fig4}}
\end{figure}
Figure \ref{shourick_contrib_Fig4} (left) show the solutions $v_4$ EMT et $i_{71}$ TS of heterogeneous DDM EMT($\Delta t=2.10^{-4}$)-TS($\Delta T=2.10^{-2}$) with comparison with the DAE solution on monodomain. We proceed to a  jump in amplitude at $t=0.04$  for the source voltage. Figure \ref{shourick_contrib_Fig4} (right) gives the $log_{10}$ of the error between two consecutive RAS iterates at $t=0.02$. It shows a  linear convergence behavior and can therefore be accelerated by the Aitken's accelerating of the convergence technique after $9$ iterates needed to numerically construct the error operator $P$.

\section{Conclusion \label{shourick_contrib_Sec5}}
A Schwarz heterogeneous DDM was used to co-simulate an RLC electrical  circuit where a part of the domain is modeled with EMT modeling and the other part with TS modeling. We showed the convergence/divergence property of the homogeneous  DDM EMT-EMT and TS-TS and of the heterogeneous DDM TS-EMT, with or without overlap and we use the pure linear divergence/convergence  of the method to accelerate it toward the searched solution with the Aitken's acceleration of the convergence  technique. The domain partitioning is only based on connectivity considerations since we want, in the long term, for the electrical network, to take advantage of the two TS and EMT representations on the overlap in order to identify the loss of information between the two models. We would like then to use this knowledge to work on other transmission conditions than Dirichlet to conserve some invariants such as electrical power.
\vspace*{-1.cm}
\bibliographystyle{plain}
\bibliography{DD26_EMTTS_v1.bbl} 
\end{document}